\documentclass{elsart}
\usepackage{amsfonts}
\usepackage{amssymb}
\usepackage{latexsym}
\usepackage{amsmath}
\usepackage{graphics}

\newtheorem{theorem}{Theorem}[section]

\newtheorem{lemma}[theorem]{Lemma}
\newtheorem{proposition}[theorem]{Proposition}

\newtheorem{open problem}[theorem]{Open Problem}

\begin{document}
\begin{frontmatter}

\title{A bound on the scrambling index of a primitive matrix using Boolean rank}

\author[akelbek,fital]{Mahmud Akelbek\corauthref{cor}}, %,Steve Kirkland}%\thanksref{NSERC}}
\ead{am44@txstate.edu}
\author[fital]{Sandra Fital},
\ead{sfitalakelbek@weber.edu}
\author[akelbek]{Jian Shen\thanksref{NSERC}}
\ead{js48@txstate.edu}

\address[akelbek]{Department of Mathematics, Texas State University, San Marcos, TX 78666
}
\address[fital]{Department of Mathematics, Weber State University, Ogden, UT 84408}
%\address[fital]{Department of Mathematics, Weber State University, Ogden, UT 84408}

\corauth[cor]{corresponding author. }
%\ead{akelbek@math.uregina.ca}}

%\author{Chun-Hua Guo\thanksref{NSERC}}
\thanks[NSERC]{Research
supported in part by NSF (CNS 0835834), Texas Higher Education Coordinating Board
(ARP 003615-0039-2007).}

%\thanks[NSERC]{This work was supported in part by a grant from the Natural Sciences and Engineering Research Council of Canada.}

\begin{abstract}  %The scrambling index of a primitive digraph $D$ is the smallest positive integer $k$
%such that for every pair of vertices $u$ and $v$, there is a  vertex $w$ such that we can get to $w$ from %$u$ and $v$ in $D$ by
%directed walks of length $k$; it is denoted  by $k(D)$.

The scrambling index of an $n\times n$ primitive matrix $A$ is the smallest positive integer $k$ such that $A^k(A^{t})^k=J$, where $A^t$ denotes the transpose of $A$ and $J$ denotes the $n\times n$ all ones matrix.
%For a primitive stochastic matrix $S$, upper bounds on the second largest modulus of an eigenvalue of $S$ are very important, because they determine the asymptotic rate of convergence of the sequence of powers of the  corresponding matrix. By using scrambling index one can provide an attainable upper bound on the second largest modulus of eigenvalues of primitive matrix.
%In this paper, we introduce the definition of the scrambling index for a primitive digraph.
 %The scrambling index of a primitive matrix $A$ can also defined as the smallest positive integer $k$ such that $A^k(A^t)^k=J$,
For an $m\times n$ Boolean matrix $M$, its {\it Boolean rank} $b(M)$ is the smallest positive integer $b$ such that $M=AB$ for some $m \times b$ Boolean matrix $A$ and $b\times n$ Boolean matrix $B$.
%
%For an $m\times n$ Boolean matrix $M$, its {\it Boolean rank} $b(M)$ is the smallest positive integer $b$ %such that for some $m \times b$ Boolean matrix $A$ and $b\times n$ Boolean matrix $B$, $M=AB$.
%The Boolean rank of the zero matrix is defined to be zero. $M=AB$ is called a {\it Boolean rank %factorization} of $M$.
In this paper, we give an upper bound on the scrambling index of an $n\times n$ primitive matrix $M$ in terms of its Boolean rank $b(M)$. Furthermore we characterize all primitive matrices that achieve the upper bound.

%In this paper, we give a bound on the scrambling index for primitive matrix $A$ in terms of its Boolean rank $b(A)$,

\end{abstract}

\noindent

\end{frontmatter}

\section {Introduction }
\noindent

For terminology and notation used here we follow \cite{RH}. A matrix $A$ is called {\it nonnegative} if all its elements are nonnegative, and denoted by $A \geq 0$. A
matrix $A$ is called {\it positive} if all its elements are positive, and denoted by $A >0$. For an $m\times n$ matrix $A$, we will denote its $(i,j)$-entry by $A_{ij}$, its $i$th row by $A_{i.}$, and its $j$th column by $A_{.j}$. For $m\times n$ matrices $A$ and $B$, we say that $B$ is dominated by $A$ if
$B_{ij}\leq A_{ij}$ for each $i$ and $j$, and denote $B\leq A$. We denote the $m\times n$ all ones matrix by $J_{m,n}$ (and by $J_{n}$ if $m=n$), The $m\times n$ all zeros matrix by $O_{m,n}$, the all ones $n$-vector by $j_n$, the $n\times n$ identity matrix by $I_n$, and its $i$th column by $e_i(n)$. The subscripts $m$ and $n$ will be omitted whenever their values are clear from the context.

For an $n \times n$ nonnegative matrix $A=(a_{ij})$, its digraph, denoted by  ${\it D(A)} $, is the digraph with
vertex set $V(D(A))=\{ 1, 2,\cdots, n\}$, and $(i,j)$ is an arc of $D(A)$ if and only if $a_{ij} \not=0$.   Then, for a positive integer $r \geq 1$,  the $(i,j)$-th entry of the matrix $A^r$ is positive if and only if  $i \overset{r}{\rightarrow} j$ in the digraph $D(A)$. Since most of the time we are only interested in the existence of such walks, not the number of different directed walks from vertex $i$ to vertex $j$, we interpret $A$ as a Boolean $(0,1)$-matrix,  unless stated otherwise. A {\it Boolean $(0,1)$-matrix} is a matrix with only $0$'s and $1$'s as its entries. Using {\it Boolean arithmetic}, $(1+1=1,\ 0+0=0,\ 1+0=1)$, we have that $AB$ and $A+B$ are Boolean $(0,1)$-matrices if $A$ and $B$ are.

Let $D=(V,E)$ denote a {\it digraph} (directed graph) with vertex set $V=V(D)$, arc set $E=E(D)$ and order
$n$. Loops are permitted but multiple arcs are not. A $u \rightarrow v$ {\it walk} in a digraph $D$ is a sequence of vertices $u,u_1,\ldots, u_t, v\ \in V(D)$ and a sequence of arcs
$(u,u_1),(u_1,u_2), \ldots, (u_t, v)\ \in\ E(D)$, where the vertices and arcs are not necessarily distinct. %A {\it closed walk} is a $u \rightarrow v$ walk where $u=v$. A {\it cycle} is a closed $u \rightarrow v$ walk with distinct vertices except for $u=v$.
%The {\it length of a walk}  $W$  is the number of arcs in $W$.
%A digraph $D$ is {\it strongly connected} if, for every ordered pair of vertices $u,\ v \in V(D)$, there is a walk from vertex $u$ to vertex $v$ in $D$.
We shall use the notation $u \rightarrow v$ and $u \nrightarrow v$ to denote, respectively, that there is an arc from vertex $u$ to vertex $v$ and that there is no such an arc. Similarly, $u\overset{k}{\rightarrow}v$ and
$u\overset{k}{\nrightarrow}v$ denote, respectively, that there is a directed walk of length $k$ from vertex $u$ to vertex $v$, and that there is no such a walk.
%The notation $ u \overset{k}{\longrightarrow} v  $ is used to indicate that there is a $u \rightarrow v$ walk of length $k$. %The {\it distance} from vertex $ u$ to vertex $v$ in $D$, is the length of a shortest walk from $u$ to $v$, and denoted by $d(u,v).$
%The {\it diameter} of a strongly connected digraph $D$ is $ \max \{d(u,v): u,v \in V(D)\}$, and denoted by $d(D)$.  A {\it $p$-cycle}  is a cycle of length $p$, denoted $C_p$. If the digraph $D$ has at least one cycle, the length of a shortest cycle in $D$ is called the {\it girth} of $D$, denoted $s(D)$.
%The number of arcs entering (leaving) a vertex $u$ is called the {\it in-degree (out-degree)} of $u$, denoted $deg^-(u)$ ($deg^+(u))$.

A digraph $D$ is called {\it primitive} if for some positive integer $t$ there is a walk of length exactly $t$ from each vertex $u$
to each
vertex $v$. If $D$ is primitive the smallest such $t$ is called the {\it exponent} of $D$, denoted by $ \exp(D)$. Equivalently, a square nonnegative matrix $A$ of order  $n$ is called {\it primitive} if there exists a positive integer $r$ such that $A^r>0$. The minimum such $r$ is called the {\it exponent} of $A$, and denoted by $\exp(A)$. Clearly $\exp(A)=\exp(D(A))$. There are numerous results on the exponent of primitive matrices \cite{RH}. %and primitive digraphs, see

%For vertices $u, v$ and $w$ of a digraph $D$, if $(u,w), (v,w) \in E(D)$, then vertex $w$ is called a {\it common out-neighbour} of vertices $u$ and $v$.
The {\it scrambling index} of a primitive digraph $D$ is the smallest positive integer $k$ such that for every pair of vertices $u$ and $v$, there exists some vertex $w=w(u,v)$ (dependent of $u$ and $v$) such that $u\overset{k}{\rightarrow} w$  and
$v\overset{k}{\rightarrow} w$ in $D$. %In other words, it is the smallest positive integer $k$ such that each pair of vertices has a common out-neighbour in $D^k$.
The scrambling index of $D$ is denoted by $k(D)$. For $u,v \in V(D)$ $(u\neq v)$, we define the {\it local scrambling index of $u$ and $v$} as
$$k_{u,v}(D)=\min\{ k:\  u \overset{k}{\rightarrow} w \mbox{ and }  v \ \overset{k}{\rightarrow} w \mbox{ for some $w\in V(D)$} \}.$$
Then
$$ k(D)=\underset{u,v\in V(D)}{\max}\{k_{u,v}(D)\}.$$

An analogous definition for scrambling index can be given for nonnegative matrices. The {\it scrambling index} of a primitive matrix $A$, denoted by $k(A)$, is the smallest positive integer $k$ such that any two rows of $A^k$ have at least one positive element in a coincident position. The scrambling index of a primitive matrix
$A$ can also be equivalently defined as the smallest positive integer $k$ such that
$A^k(A^{t})^k=J$, where $A^t$ denotes the transpose of $A$. If $A$ is the adjacency matrix of a primitive digraph $D$, then $k(D)=k(A)$. As a result, throughout the paper, where no confusion occurs, we use the digraph $D$ and the adjacency matrix $A(D)$ interchangeably.

%by $\it J, \ O$, and $I$ the  {\it all $1$'s matrix}, the  {\it all $0$'s matrix} and the {\it identity matrix}, respectively. T

In \cite{AKEL1} and \cite{AKEL2}, Akelbek and Kirkland obtained an upper bound on the scrambling index of a primitive digraph $D$ in terms of the order and girth of $D$, and gave a characterization of the primitive digraphs with the largest scrambling index.

\begin{theorem}{\bf \cite{AKEL1}} %\label{MAIN}
Let $D$ be a primitive digraph with $n$ vertices and girth $s$. Then
$$k(D) \leq n-s+\left\{ \begin{array}{cc}
                               (\frac{s-1}{2})n, &\mbox{when $s$ is odd,}\\
                               (\frac{n-1}{2})s, &\mbox{when $s$ is even.}
                              \end{array}\right.
$$
\end{theorem}

When $s =n-1$, an upper bound on $k(D)$ in terms of the order of a primitive digraph $D$ can be achieved \cite{AKEL1}.  We state the theorem in terms of primitive matrices below.
\begin{theorem}{\bf \cite{AKEL1}}\label{pro12} Let $A$ be a primitive matrix of order $n\geq 2$. Then
\begin{equation}\label{nnnn}
 k(A) \leq \left\lceil \frac{(n-1)^2+1}{2}\right\rceil.
\end{equation}
Equality holds in (\ref{nnnn}) if and only if there is a permutation matrix $P$ such that $PAP^t$ is one of the following matrices
\[  \begin{array}{cccc}
    W_2= \left[\begin{array}{cc}
              1 & 1 \\
              1 & 0
          \end{array} \right] & \mbox{ or } &
                J_2=\left[ \begin{array}{cc}
              1 & 1 \\
              1 & 1
          \end{array} \right] & \mbox{ when $n=2$},
\end{array}
\]
$$ W_n= \left[ \begin{array}{cccccc}
              0        & 1      & 0      &\cdots    &  0       &   0      \\
              0        & 0      & 1      & 0        & \cdots   &   0      \\
              \vdots   & \vdots & \ddots & \ddots   & \cdots   & \vdots    \\
              0        & 0      & \cdots & 0        & 1        &   0       \\
              1        & 0      & \cdots & \cdots   & 0        &   1       \\
              1        & 0      & \cdots & \cdots   & 0        &   0

          \end{array} \right]   \mbox{when $ n \geq 3$.}
     $$
\end {theorem}

The digraph $D(W_n)$ is called the Wielandt graph and denoted by $D_{n-1,n}$. It is a digraph with a Hamilton cycle $1\rightarrow 2\rightarrow\cdots \rightarrow n \rightarrow 1$ together with an arc from vertex $n-1$ to vertex $1$.
%Let $D=D(W_n)$, then $D=D_{n-1,n}$, and we label $D_{n-1,n}$ as in Figure \ref{ }. % In \cite{AKEL1} we gave a detailed description for $D_{s,n}$, where $2\leq s\leq n$. Now we are considering a special case that $s=n-1$ and we
%\begin{figure}[h!]
%\begin{center}
%\scalebox{0.5}{\includegraphics{za1.eps}}
%\caption{ $D_{s,n}$ }
%\label{dsnwwaa}
%\end{center}
%\end{figure}
For simplicity, let $h_n=\left\lceil \frac{(n-1)^2+1}{2}\right\rceil$. The next proposition gives some information about the Wielandt graph $D_{n-1,n}$. %Here are some results on the scrambling index of Wielandt graph $D_{n-1,n}$ \cite{AKEL1}.

\noindent
\begin{proposition}{\bf \cite{AKEL1}}\label{prop23}
For $D_{n-1,n}$, where $n\geq 3$,
\begin{description}
  \item[$(a)$] $k_{n,\lfloor\frac{n}{2}\rfloor}(D_{n-1,n})=h_n$, and for all other pairs of vertices $u$ and $v$ of $D_{n-1,n}$, $k_{u,v}(D_{n-1,n})<h_n$.
  \item[$(b)$] There are directed walks from vertices $n$ and $\lfloor\frac{n}{2}\rfloor$ to vertex $1$ of length $h_n$, that is $n\overset{h_n}{\rightarrow} 1$ and $\lfloor\frac{n}{2}\rfloor \overset{h_n}{\rightarrow} 1$.
\end{description}
  \end{proposition}

%\begin{theorem} Let $D$ be a primitive digraph of order $n\geq 2$. Then
%\begin{equation}\label{nnnn}
% k(D) \leq \left\lceil \frac{(n-1)^2+1}{2}\right\rceil.
%\end{equation}
%Equality holds if and only if  $D=D_{n-1,n}$.
%\end {theorem}
%In \cite{GREG}, D.~Gregory et al. gave the upper bound on the exponent in terms of the Boolean rank of a primitive matrix $A$.
For an $m\times n$ Boolean matrix $M$, we define its {\it Boolean rank} $b(M)$ to be the smallest positive integer $b$ such that for some $m \times b$ Boolean matrix $A$ and $b\times n$ Boolean matrix $B$, $M=AB$.
The Boolean rank of the zero matrix is defined to be zero. $M=AB$ is called a {\it Boolean rank factorization} of $M$.

In \cite{GREG}, Gregory, Kirkland and Pullman obtained an upper bound on the exponent of primitive Boolean matrix in terms of Boolean rank.

\begin{proposition}{\bf \cite{GREG}} Suppose that $n\geq 2$ and that $M$ is an $n\times n$ primitive Boolean matrix with $b(M)=b$. Then
\begin{equation}\label{greg1}
\exp(M) \leq (b-1)^2+2.
\end{equation}
\end{proposition}

In \cite{GREG}, Gregory, Kirkland and Pullman also gave a characterization of the matrices for which equality holds in (\ref{greg1}). In \cite{LIU1}, Liu, You and Yu gave a characterization of primitive matrices $M$ with Boolean rank $b$  such that $\exp(M)=(b-1)^2+1$.

In this paper, we give an upper bound on the scrambling index  of a primitive matrix $M$ using  Boolean rank $b=b(M)$, and characterize all Boolean primitive matrices that achieve the upper bound.

\section{Main Results}

%\subsection{A bound on the scrambling index}

We start with a basic result.

\begin{lemma}\label{pro11}
 Suppose that $A$ and $B$ are $n\times m$ and $ m\times n$ Boolean matrices respectively, and that neither has a zero line. Then

 $(a)$ $AB$ is primitive if and only if $BA$ is primitive.

 $(b)$ If $AB$ and $BA$ are primitive, then
\begin{equation}
|k(AB)-k(BA)|\leq 1.
\end{equation}
\end{lemma}

{\bf Proof.} Part (a) was proved by Shao \cite{SHAO1}. We only need to show part $(b)$. Since $AB$ and $BA$ are primitive matrices, $A$ and $B$ has no zero rows. Then $AA^t\geq I_n$ and $BJ_nB^t=J_m$. Suppose $k(AB)=k$. By the definition of scrambling index
$$ (AB)^k((AB)^t)^k=J_n.$$
Then
\begin{eqnarray*}
(BA)^{k+1}((BA)^t)^{k+1} &=& B(AB)^kAA^t((AB)^t)^kB^t\geq B(AB)^kI_n((AB)^t)^kB^t \\
          & =  & B(AB)^k((AB)^t)^kB^t=BJ_nB^t=J_m.
\end{eqnarray*}

Thus $k(BA)\leq k+1=k(AB)+1$. The result follows by exchanging the roles of $A$ and $B$. $\Box$

\begin{proposition}{\bf \cite{LIU1}}\label{zeko1}
Let $M$ be an $n\times n$ primitive Boolean matrix, and $M=AB$ be a Boolean rank factorization of $M$. Then neither $A$ nor $B$ has a zero line.
\end{proposition}

%\begin{proposition}\label{zeko1}
%Let $M$ be an $n\times n$ primitive Boolean matrix, and $M=AB$ is a Boolean rank factorization of $M$. Then %neither $A$ nor $B$ has a zero line.
%\end{proposition}

%{\bf Proof.} Suppose the Boolean rank of $M$ is $b$, and $M=AB$ is a Boolean rank factorization with $A$ and $B$ are $m\times b$ and $b\times n$ matrices, respectively. If matrix $A$ has a zero row, then $M$ has a zero row, and $M$ is not primitive. A contradiction. If matrix $A$ has zero a column, then the Boolean rank of $M$ is smaller than $b$. Again a contradiction. We can apply the same argument to matrix $B$. Therefore neither $A$ nor $B$ has a zero line.  $\Box$

 %, and the second one

\begin{theorem}\label{ggg}
Let $M$ be an $n\times n$ $(n\geq 2)$ primitive matrix with Boolean rank $b(M)=b$. Then
\begin{equation}\label{scram1}
k(M)\leq \left\lceil\frac{(b-1)^2+1}{2}\right\rceil +1.
\end{equation}
\end{theorem}

{\bf Proof.} Let $M=AB$ be a Boolean rank factorization of $M$, where $A$ and $B$ are $n\times b$ and $b\times n$ Boolean matrices respectively. Then by Lemma \ref{zeko1} neither $A$ nor $B$ has a zero line. By lemma \ref{pro11}, we have
$$k(M)=k(AB)\leq k(BA)+1.$$
Since $BA$ is primitive and $BA$ is a $b\times b$ matrix, by Theorem \ref{pro12},
$$k(BA)\leq \left\lceil\frac{(b-1)^2+1}{2}\right\rceil,$$
from which Theorem \ref{ggg} follows. $\Box$

From (\ref{nnnn}) we see that no matrix of full Boolean rank $n$ can attain the upper bound in (\ref{scram1}). Further, since the only $n\times n$ primitive Boolean matrix with Boolean rank $1$ is $J_n$, no matrix of Boolean rank $1$ can attain the upper bound in (\ref{scram1}). Thus we may assume that $2\leq b\leq n-1$.

For simplicity, let
$$h=\left\lceil\frac{(b-1)^2+1}{2}\right\rceil.$$

%\subsection{The case when $b\geq 3$}

Recall from Theorem \ref{pro12} that $k(W_b)=h$. We first make some observations about $W_b$. Recall that $D=D(W_b)$ is the Wielandt graph $D_{b-1,b}$ with $b$ vertices.

\begin{lemma}\label{scram2}
If $b\geq 3$, then the zero entries of $(W_b)^{h-1}(W_b^t)^{h-1}$ occur only in the
$(b,\lfloor \frac{b}{2}\rfloor)$ and $(\lfloor \frac{b}{2}\rfloor, b)$ positions.
\end{lemma}

{\bf Proof.} By Proposition \ref{prop23} we know that $k_{b,\lfloor\frac{b}{2}\rfloor}(D_{b-1,b})=h$, and for all other pairs of vertices $u$ and $v$, $k_{u,v}(D_{b-1,b})<h$. Therefore in $W_b^{h-1}$ every pair of rows intersect with each other except rows $b$ and $\lfloor\frac{b}{2}\rfloor$.  Thus the only zero entries of $(W_b)^{h-1}(W_b^t)^{h-1}$  are in the
$(b,\lfloor \frac{b}{2}\rfloor)$ and $(\lfloor \frac{b}{2}\rfloor, b)$ positions. $\Box$

%$$ (BA)^{k+1}((BA)^t)^){k+1}=B(AB)^kAA^t((AB)^t)^)kB^t\geqB(AB)^kI_n((AB)^t)^)kB^t$$

For an $n\times n$ $(n\geq2)$ matrix $A$, let $A(\{i_1,i_2\}, \{j_1,j_2\})$ be the submatrix of $A$ that lies in the rows $i_1$ and $i_2$ and the columns $j_1$ and $j_2$.

\begin{lemma}\label{scram23}
For $b\geq 3$, $W_b^{h-1}(\{\lfloor\frac{b}{2}\rfloor,b\}, \{b-1,b\})$ is either
$\left[ \begin{array}{cc}
               1 & 0\\
               0 & 1
            \end{array} \right]$ or
$\left[ \begin{array}{cc}
               0 & 1\\
               1 & 0
            \end{array} \right]$.
\end{lemma}

{\bf Proof.}  By Proposition \ref{prop23}, we know that $k_{\lfloor\frac{b}{2}\rfloor, b}(D_{b-1,b})=h$ and
$\lfloor\frac{b}{2}\rfloor \overset{h}{\rightarrow} 1$  and $b\overset{h}{\rightarrow} 1.$
From the digraph $D_{b-1,b}$, we know that the directed walks of length $h$ from vertices $\lfloor\frac{b}{2}\rfloor$ and $b$ to vertex $1$ is either
$$  \lfloor\frac{b}{2}\rfloor \overset{h-1}{\rightarrow}b-1 \overset{1}{\rightarrow} 1, $$
       $$b\overset{h-1}{\rightarrow} b \overset{1}{\rightarrow} 1,
$$ or
$$  \lfloor\frac{b}{2}\rfloor \overset{h-1}{\rightarrow}b \overset{1}{\rightarrow} 1, $$
       $$b\overset{h-1}{\rightarrow} b-1 \overset{1}{\rightarrow} 1.
$$

For the first case, if $\lfloor\frac{b}{2}\rfloor \overset{h-1}{\rightarrow}b-1$ and $b\overset{h-1}{\rightarrow} b$, then $b\overset{h-1}{\nrightarrow}b-1$ and $\lfloor\frac{b}{2}\rfloor \overset{h-1}{\nrightarrow} b$. Otherwise it contradicts to $k_{\lfloor\frac{b}{2}\rfloor,b}(D_{b-1,b})=h$. Similarly, for the second case if $\lfloor\frac{b}{2}\rfloor \overset{h-1}{\rightarrow}b$ and $b\overset{h-1}{\rightarrow} b-1$, then $b \overset{h-1}{\nrightarrow}b$ and $\lfloor\frac{b}{2}\rfloor\overset{h-1}{\nrightarrow} b-1$. The result follows by applying these to the matrix $W_{b}^{h-1}$.
$\Box$

%Let $u=e_{b-1}(b)+e_{b}(b)$ and $u'=e_{\lfloor\frac{b}{2}\rfloor}(b)+e_b(b).$

%\begin{lemma}\label{lemw2} For $b\geq 3$, $ (W_b)^{h-1}(I+uu^t)((W_b^t)^{h-1}=J.$

%\end{lemma}

%{\bf Proof.}

\begin{theorem}\label{theo12}
Suppose $M$ is an $n\times n$ primitive Boolean matrix with $3\leq b=b(M)\leq n-1$. Then
$$k(M)=\left\lceil\frac{(b-1)^2+1}{2}\right\rceil+1$$
if and only if $M$ has a boolean rank factorization $M=AB$, where $A$ and $B$ have the following properties:
\begin{description}
  \item    \mbox{\rm{(i)}} $BA=W_b$,
  \item  \mbox{\rm{(ii)}} some row of $A$ is $e_{\lfloor\frac{b}{2}\rfloor}^t(b)$, some row of $A$ is $e_b^t(b)$, and
  \item  \mbox{\rm{(iii)}} no column of $B$ is $e_{b-1}(b)+e_{b}(b)$.
 \end{description}
\end{theorem}

{\bf Proof.} First suppose $M$ is primitive with $k(M)=h+1$, and $M=\tilde{A} \tilde{B}$ is a Boolean rank factorization of $M$. By Lemma \ref{pro11}, $\tilde{B}\tilde{A}$ is primitive and $k(\tilde{B}\tilde{A})\geq h$. But $\tilde{B}\tilde{A}$ is a $b\times b$ matrix. By Theorem \ref{pro12}, $k(\tilde{B}\tilde{A})\leq h$. Therefore $k(\tilde{B}\tilde{A})=h$. Also by Theorem \ref{pro12}, there is a permutation matrix $P$ such that $P\tilde{B}\tilde{A}P^t=W_b$. Let $B=P\tilde{B}$ and $A=\tilde{A}P^t$. Then $AB=\tilde{A}P^tP\tilde{B}=\tilde{A}\tilde{B}=M$. Thus $A$ and $B$ satisfy condition $\mbox{(i)}$.

Since $M$ is primitive, we have $\sum_{i=1}^b A_{.i}=j_n=\sum_{i=1}^b B_{i.}^t$. Since $k(M)=h+1$, the matrix $M^h$ must have two rows that do not intersect. Without lost of generality, suppose rows $p$ and $q$ of $M^h$ do not intersect. Then entries in the $(p,q)$ and $(q,p)$ positions of $M^h(M^t)^h$ are zero. Since matrix $B$ has no zero row, we have $BB^t\geq I_b$. Thus
\vspace{-.2cm}
$$\begin{array}{ll}
    & M^h(M^t)^h    \\
=& (AB)^h((AB)^t)^h=A(BA)^{h-1}BB^t((BA)^t)^{h-1}A^t \\
           =  &A(W_b)^{h-1}BB^t(W_b^t)^{h-1}A^t\\
          \geq & A(W_b)^{h-1}I_b(W_b^t)^{h-1}A^t=A(W_b)^{h-1}(W_b^t)^{h-1}A^t \\
          =& AZA^t\\
          =&\left[ J_{n,\lfloor\frac{b}{2}\rfloor-1}\vline \displaystyle \sum_{i=1}^{b-1}A_{.i}\vline J_{n,b-\lfloor\frac{   b }{2}\rfloor-1}\vline \displaystyle \sum^b_{\underset{{i\neq\lfloor\frac{b}{2}\rfloor}}{i=1}}A_{.i} \right]  A^t\\
          =& j_n \left( \displaystyle \sum_{i=1}^{\lfloor\frac{b}{2}\rfloor-1} A_{.i} \right)^t+\left(\displaystyle \sum_{i=1}^{b-1} A_{   .i }\right)(A_{.\lfloor\frac{b}{2}\rfloor})^t +j_n \left(\displaystyle \sum_{i=\lfloor\frac{b}{2}\rfloor+1}^{b-1} A_{.i}  \right)^t+\left(\displaystyle \sum^b_{\underset{{i\neq\lfloor\frac{b}{2}\rfloor}}{i=1}} A_{.i}\right) (A_{.b})^t,
  \end{array}$$
%%%%  This one is the old one
%  \noindent
% \begin{eqnarray*}
%  \noindent
%    M^h(M^t)^h    &=& (AB)^h((AB)^t)^h=A(BA)^{h-1}BB^t((BA)^t)^{h-1}A^t \\
%          & =  &A(W_b)^{h-1}BB^t(W_b^t)^{h-1}A^t\\
%          & \geq & A(W_b)^{h-1}I_b(W_b^t)^{h-1}A^t=A(W_b)^{h-1}(W_b^t)^{h-1}A^t \\
%          &=& AZA^t\\
%          &=&\left[ J_{n,\lfloor\frac{b}{2}\rfloor-1}\vline \sum_{i=1}^{b-1}A_{.i}\vline J_{n,b-\lfloor\frac{   b }{2}\rfloor-1}\vline\sum^b_{\underset{{i\neq\lfloor\frac{b}{2}\rfloor}}{i=1}}A_{.i} \right]  A^t\\
%          &=& j_n \left(\sum_{i=1}^{\lfloor\frac{b}{2}\rfloor-1} A_{.i} \right)^t+\left(\sum_{i=1}^{b-1} A_{   .i }\right)(A_{.\lfloor\frac{b}{2}\rfloor})^t +j_n \left(\sum_{i=\lfloor\frac{b}{2}\rfloor+1}^{b-1} A_{.i}  \right)^t+\left(\sum^b_{\underset{{i\neq\lfloor\frac{b}{2}\rfloor}}{i=1}} A_{.i}\right) (A_{.b})^t,
%  \end{eqnarray*}
where $Z=(W_b)^{h-1}(W_b^t)^{h-1}$ is the $b\times b$ matrix which has zero entries only in the $(\lfloor\frac{b}{2}\rfloor, b)$ and
$(b, \lfloor\frac{b}{2}\rfloor)$ positions. Since $AZA^t$ is dominated by $M^h(M^t)^h$ and $M^h(M^t)^h$ has zero entries in the $(p,q)$ and $(q,p)$ positions, the entries in the $(p,q)$ and $(q,p)$ positions of $AZA^t$ are also zero. Thus
\begin{equation}\label{pq1}
\sum_{i=1}^{\lfloor\frac{b}{2}\rfloor-1}A_{qi}+\left(\sum_{i=1}^{b-1}A_{pi}\right)A_{q\lfloor\frac{b}{2}\rfloor}+\sum_{i=\lfloor\frac{b}{2}\rfloor+1}^{b-1}A_{qi}+\left(\sum_{\underset{{i\neq\lfloor\frac{b}{2}\rfloor}}{i=1}}^{b}A_{pi}\right)A_{qb}=0
\end{equation}
and
\begin{equation}\label{pq2}
 \sum_{i=1}^{\lfloor\frac{b}{2}\rfloor-1}A_{pi}+\left(\sum_{i=1}^{b-1}A_{qi}\right)A_{p\lfloor\frac{b}{2}\rfloor}+\sum_{i=\lfloor\frac{b}{2}\rfloor+1}^{b-1}A_{pi}+\left(\sum_{\underset{{i\neq\lfloor\frac{b}{2}\rfloor}}{i=1}}^{b}A_{qi}\right)A_{pb}=0.
\end{equation}

Then $A_{qi}=0$ and $A_{pi}=0$ for $i=1,\cdots, b-1$ and $i\neq \lfloor\frac{b}{2}\rfloor$. Substitute these back to (\ref{pq1}) and (\ref{pq2}), we have
\begin{equation}\label{pq3}
A_{q\lfloor\frac{b}{2}\rfloor}A_{p\lfloor\frac{b}{2}\rfloor}+A_{qb}A_{pb}=0.
\end{equation}
%and
%\begin{equation}\label{pq4}
%A_{p\lfloor\frac{b}{2}\rfloor}A_{q\lfloor\frac{b}{2}\rfloor}+A_{pb}A_{qb}=0,
%\end{equation}

If $A_{q\lfloor\frac{b}{2}\rfloor}\neq 0$, then $A_{p\lfloor\frac{b}{2}\rfloor}=0$. Since every row of $A$ is nonzero, we have $A_{pb}\neq 0$. By (\ref{pq3}), $A_{qp}=0$. Therefore some rows of $A$ is $e_{\lfloor\frac{b}{2}\rfloor}^t(b)$ and some row of $A$ is $e_{b}^t(b)$. This concludes $\mbox{(ii)}$.

We claim $B$ can not have a column which is equal to $u$. Otherwise, suppose some column of $B$ is $u$.
Since $B$ has no zero row, by Proposition \ref{zeko1}, $BB^t\geq I_b+uu^t$. Thus
%Since $BA=W_b$ by $\mbox{(i)$, and $A$ has no zero row, each column of $B$ is dominated by a column of $W_b$. %Similarly, each row of $A$ is dominated by a row of $W_b$. Thus each column of $B$ is in the set
%$S_1=\{e_1(b), e_2(b), \cdots, e_b(b), u\}$, where $u=e_{b-1}(b)+e_b(b)$, and each row of $A$ is in the set
%$S_2=\{e_1^t(b), e_2^t(b), \cdots, e_b^t(b), v^t\}$, where $v=e_1(b)+e_b(b)$.
%Next, we note that for each $1\leq i\leq b$, the outer product $B_{.i}A_{i.}$ is dominated by $W_b$. Since each such $B_{.i}$ and $A_{i.}$ must be in $S_1$ and $S_2$ respectively, we find that $(B_{.i},A_{i.})$ must be one of the following pairs: $(e_i,e_{i+1}^t)$, $1\leq i\leq b-1$, $(e_{b-1}, e_{1}^t)$, $(e_b, e_1)$, $(u,e_1^t)$, and $(e_{b-1}, v^t)$.
%Since for each $1\leq i \leq b-1$, some outer product $B_{.i}A_{i.}$ have a $1$ in the $(i,i+1)$ positions, hence for some $k_i$, $(B_{k_i},A_{k_i})=(e_i,e_{i+1}^t)$. Some outer product $B_{.i}A_{i.}$ have $1$'s in the $(b,1)$ and $(b-1,1)$ positions, hence for some $k_j$, $(B_{k_j},A_{k_j})$ is one of
%$(e_{b-1},e_{1}^t)$, $(e_b,e_{1}^t)$, $(e_{b-1}, v^t)$ or $(u,e_1^t)$.
%We claim $B$ can not have a column which is equal to $u$. Otherwise, suppose some column of $B$ is $u$, then
%each $e_1, e_2, \cdots, e_{b-1}$ and $u$ is a column of $B$, and
%$$BB^t=\sum_{i=1}^{b-1} e_ie_i^t+uu^t=I+uu^t.$$
\begin{eqnarray*}
\noindent
M^h(M^t)^h &=& (AB)^h((AB)^t)^h=A(BA)^{h-1}BB^t((BA)^t)^{h-1}A^t \\
          & = &A(W_b)^{h-1}BB^t(W_b^t)^{h-1}A^t\\
          & \geq & A(W_b)^{h-1}(I_b+uu^t)(W_b^t)^{h-1}A^t\\
          &=& A[(W_b)^{h-1}(W_b^t)^{h-1}+(W_b^{h-1}u)(W_b^{h-1}u)^t]A^t.
                    \end{eqnarray*}
By lemma \ref{scram2}, $W_b^{h-1}(\{\lfloor\frac{b}{2}\rfloor,b\}, \{b-1,b\})$ is either
$\left[ \begin{array}{cc}
               1 & 0\\
               0 & 1
            \end{array} \right]$ or
$\left[ \begin{array}{cc}
               0 & 1\\
               1 & 0
            \end{array} \right]$. Then $W_b^{h-1}u \geq e_{\lfloor\frac{b}{2}\rfloor}(b)+e_b(b)$. By Lemma \ref{scram2}, the zero entries of $W_b^{h-1}(W_b^t)^{h-1}$ are in the $(b,\lfloor \frac{b}{2}\rfloor)$ and $(\lfloor \frac{b}{2}\rfloor, b)$ positions. Therefore $ W_b^{h-1}(W_b^t)^{h-1}+(W_b^{h-1}u)(W_b^{h-1}u)^t=J_b$. Since $A$ has no zero lines, we have $ M^h(M^t)^h=AJ_bA^t=J_n$, which is a contradiction to $k(M)=h+1$. This proves $\mbox{(iii)}$.

Finally, suppose that $M=AB$ is a Boolean rank factorization of $M$ and $A$ and $B$ satisfy $\mbox{(i)}$, $\mbox{(ii)}$ and $\mbox{(iii)}$. By Lemma \ref{pro11}$(a)$ and Theorem \ref{pro12}, the matrix $M$ is primitive and $k(M)\leq h+1$ by Lemma \ref{pro11}$(b)$ and . But it follows from Lemma \ref{scram2} and conditions $\mbox{(i)}$, $\mbox{(ii)}$ and $\mbox{(iii)}$ that $M^h$ has zero entries. So we conclude that $k(D)=h+1$.   $\Box$

Next we will reinterpret conditions $\mbox{(i)}$, $\mbox{(ii)}$ and $\mbox{(iii)}$ of Theorem \ref{theo12} to show that if $k(M)=h+1$, then $M$ is one of the three basic types of matrices in Theorem \ref{char1}.

$$ Table\  1 \ \ ( b \geq 3) $$
\begin{tabular}{ll} \hline \\
 $M_1=
\left[ \begin{array}{cccccc|c}
              0        & J     & 0      &0       & \cdots  &0       &0       \\
             \vdots   & \vdots & \ddots &\ddots  & \ddots  &\vdots   &\vdots  \\
             0        & 0      & \cdots &J       & 0       &0        &0       \\
             0        & 0      & \cdots &0       & J       &0        &J       \\
             0        & 0      & \cdots &0       & 0       &J        &0        \\
             J        & 0      & \cdots &0       & 0       &0        &0       \\ \hline
             J        & 0      & \cdots &0       & 0       &0        &0
            \end{array} \right]$ &
$M_2=
\left[ \begin{array}{cccccc|c}
              0        & J     & 0      &0       & \cdots  &0       &0       \\
             \vdots   & \vdots & \ddots &\ddots  & \ddots  &\vdots   &\vdots  \\
             0        & 0      & \cdots &J       & 0       &0        &0       \\
             0        & 0      & \cdots &0       & J       &0        &J       \\
             0        & 0      & \cdots &0       & 0       &J        &0        \\
             J        & 0      & \cdots &0       & 0       &0        &0       \\ \hline
             J        & 0      & \cdots &0       & 0       &J        &0
            \end{array} \right]$  \\  \\
$M_3=
\left[ \begin{array}{cccccc|cc}
              0        & J     & 0      &0       & \cdots  &0       &0       & 0     \\
             \vdots   & \vdots & \ddots &\ddots  & \ddots  &\vdots   &\vdots & \vdots \\
             0        & 0      & \cdots &J       & 0       &0        &0      & 0     \\
             0        & 0      & \cdots &0       & J       &0        &J      & J     \\
             0        & 0      & \cdots &0       & 0       &J        &0      & 0    \\
             J        & 0      & \cdots &0       & 0       &0        &0      & 0    \\ \hline
             J        & 0      & \cdots &0       & 0       &0        &0      & 0   \\
             J        & 0      & \cdots &0       & 0       &J        &0      & 0
            \end{array} \right]$ &  \\ \\
\hline
\end{tabular}

\begin{theorem}\label{char1} Suppose $M$ is an $n \times n$ Boolean matrix with $b(M)=b$, where $3\leq b\leq n-1$. Then $M$ is primitive with $k(M)=h+1$ if and only if there is a permutation matrix $P$ such that $PMP^t$ has one of the forms in Table $1$. %if $3\leq b\le n-1$ or Table $2$ if $b=2$.
\end{theorem}

  In Table $1$ %and Table $2$
the rows and columns of $M_1$, $M_2$ and $M_3$ are partitioned conformally, so that each diagonal block is square, and the top left hand submatrix common to each has $b$ blocks in its partitioning.

{\bf Proof.} Suppose $M$ is primitive, $b\geq 3$, and $k(M)=h+1$. Then by Theorem \ref{theo12}$\mbox{(i)}$, $M$ has  a Boolean rank factorization $M=AB$ such that $BA=W_b$. Since $A$ has no zero row, each column of $B$ is dominated by a column of $W_b$. Similarly, each row of $A$ is dominated by a row of $W_b$. Thus each column of $B$ is in the set
$S_1=\{e_1(b), e_2(b), \cdots, e_b(b), u\}$, where $u=e_{b-1}(b)+e_b(b)$. Similarly, each row of $A$ is in the set
$S_2=\{e_1^t(b), e_2^t(b), \cdots, e_b^t(b), v^t\}$, where $v=e_1(b)+e_b(b)$. But by Theorem
\ref{theo12}$\mbox{(iii)}$, no column of $B$ is $u$. Hence each column of $B$ is in the set of $S_1'=\{e_1(b), e_2(b), \cdots, e_b(b)\}$.

Next, we note that for each $1\leq i\leq b$, the product $B_{.i}A_{i.}$ is dominated by $W_b$. Since each $B_{.i}$ and $A_{i.}$ must be in $S'_1$ and $S_2$ respectively and $(B_{.i}, A_{i.})$ must be one of the following pairs: $(e_i,e_{i+1}^t)$, $1\leq i\leq b-1$, $(e_{b-1}, e_{1}^t)$, $(e_b, e_1)$, or $(e_{b-1}, v^t)$, where $e_i=e_i(b)$ for any $i\in \{1, 2, \cdots, b\}$. Thus, for each $i$, $1\leq i \leq b-1$, $(e_i,e_{i+1}^t)=(B_{.k_i},A_{k_i.})$ for some $k_i$. Some outer product $B_{.j}A_{j.}$ has a $1$ in the $(b,1)$ position, hence $(B_{.k_b},A_{k_b.})=(e_b,e_{1}^t)$ for some $k_b$. Finally some outer product $B_{.j}A_{j.}$ must have a $1$ in the $(b-1,1)$ position, hence for some $k_{b+1}$, $(B_{.k_{b+1}},A_{k_{k+1}.})$ is one of
$(e_{b-1},e_{1}^t)$ or $(e_{b-1}, v^t)$.
%We claim that $B$ can not have a column which is equal to $u$. Consider $M^h=(AB)^h=A(W_b)^{h-1}B$. By Proposition \ref{prop23}, in
%$(W_b)^{h-1}$, only row $\lfloor\frac{b}{2}\rfloor$ and row $b$ do not intersect and all other pairs of rows  intersect with each other. Since each row of $A$ is in the set of $S_2$, then multiplying $(W_b)^{h-1}$ from right side by $A$, only interchange rows of $(W_b)^{h-1}$and sum the corresponding entries of first and last rows of $(W_b)^{h-1}$. But since some row of $A$ is $e^t_{\lfloor\frac{b}{2}\rfloor}$, and some row of $A$ is $e^t_b$, then  there are still two rows of $A(W_b)^{h-1}$ that do not intersect. By Lemma \ref{scram2}, we know that $W_b^{h-1}(\{\lfloor\frac{b}{2}\rfloor,b\}, \{b-1,b\})$ is either
%$$\left[ \begin{array}{cc}
%               1 & 0\\
%               0 & 1
%            \end{array} \right], \mbox{\  \ or\ \ }
%\left[ \begin{array}{cc}
%               1 & 0\\
%               0 & 1
%            \end{array} \right].$$
%If some column of $B$ is $u$, say column $p$, then $((W_{b})^{h-1}B)_{\lfloor\frac{b}{2}\rfloor p}=1$ and
%$((W_{b})^{h-1}B)_{bp}=1$. Thus row $\lfloor\frac{b}{2}\rfloor$ and row $b$ in $(W_{b})^{h-1}B$ intersect. This concludes every pair of rows of $A(W_{b})^{h-1}B$ intersect with each other, and $k(M)\leq h$, which is contradiction to $k(M)=h+1$. Therefore each column of $B$ is in the set of $S'_1=\{e_1(b), e_2(b), \cdots, e_b(b)\}$.
It follows from the above argument that there is an $n\times n$ permutation matrix $Q$ such that $$ BQ^t=\left[\bar{B}|\tilde{B}\right] \mbox{\ \ \ \ and \ \ \ }
                  QA=\left[\begin{array}{c}
                             \bar{A}\\ \hline
                             \tilde{A}
                            \end{array}\right], $$
where $$\bar{B}=\left[e_1j_{n_1}^t | e_2j_{n_2}^t |\cdots| e_bj_{n_b}^t\right] \mbox{\ \ and \ \ }
               \bar{A}=\left[\begin{array}{c}
                           j_{n_1}e_2^t\\ \hline
                           j_{n_2}e_3^t\\ \hline
                           \cdots \\ \hline
                           j_{n_{b-1}}e_b^t\\ \hline
                           j_{n_b}e_1^t
                        \end{array}\right]$$
for some $n_1,\cdots, n_b\geq 1$, and where each $(\tilde{B}_{.i}, \tilde{A}_{i.})$ is one of
$(e_{b-1}, e_1^t)$ or $(e_{b-1}, v^t)$. Thus $\tilde{B}$ and $\tilde{A}$ can be one of the following pairs of matrices:
\begin{eqnarray*}
     \tilde{B}_1 & =& e_{b-1}j_{m_1}^t, \   \tilde{A}_1=j_{m_1}e_1^t \mbox{\ \ for some $m_1\geq 1$};\\
     \tilde{B}_2 & =& e_{b-1}j_{m_2}^t,  \  \tilde{A}_2=j_{m_2}v^t \mbox{\ \ for some $m_2\geq 1$};\\
     \tilde{B}_3 & =&\left[ e_{b-1}j_{m_3}^t | e_{b-1}j_{p_3}^t\right],  \
              \tilde{A}_3=\left[ \begin{array}{c}
                            j_{m_3}e_{1}^t \\ \hline
                            j_{p_3}v^t
                           \end{array}\right] \mbox{\ \ for some $m_3, p_3\geq 1$}.
\end{eqnarray*}
It is now readily verified that
$$\left[ \begin{array}{c}
             \bar{A}\\ \hline
             \tilde{A_i}
            \end{array} \right]  \left[ \bar{B}| \tilde{B_i}\right]=M_i \ \ \mbox{ for $1\leq i\leq3$,}$$
so that $QMQ^t$ is one of the matrices in Table $1$.

Finally, since the Boolean rank factorization
$$M_i=\left[ \begin{array}{c}
             \bar{A}\\ \hline
             \tilde{A_i}
            \end{array} \right]  \left[ \bar{B}| \tilde{B_i}\right]$$
satisfies conditions $\mbox{(i)}$, $\mbox{(ii)}$ and $\mbox{(iii)}$ of Theorem \ref{theo12}, each $M_i$ is primitive and $k(M)=h+1$.    $\Box$

%\subsection{Case when $b=2$}

When $b(M)=2$, we have the following result.

\begin{theorem}\label{theo15}
Suppose $M$ is an $n\times n$ primitive Boolean matrix with $b(M)=b=2$. Then
$k(M)=2$ if and only if $M$ has a boolean rank factorization $M=AB$, where $A$ and $B$ have the following properties:
\begin{description}
  \item[]  \mbox{\rm{(i)}}   $BA=W_2$ or $BA=J_2$,\\
   \item[] \mbox{\rm{(ii)}} some row of $A$ is $e_{1}^t(2)$, some row of $A$ is $e_2^t(2)$, and\\
    \item[] \mbox{\rm{(iii)}} no column of $B$ is $e_{1}(2)+e_{2}(2)$.
\end{description}
 \end{theorem}

{\bf Proof.}  First suppose $M$ is primitive with $k(M)=2$, and $M=\tilde{A} \tilde{B}$ is a Boolean rank factorization of $M$. By Lemma \ref{pro11}, $\tilde{B}\tilde{A}$ is primitive and $k(\tilde{B}\tilde{A})\geq 1$. But $\tilde{B}\tilde{A}$ is a $2\times 2$ matrix. By Theorem \ref{pro12}, $k(\tilde{B}\tilde{A})\leq 1$. Therefore $k(\tilde{B}\tilde{A})=1$. Also by Theorem \ref{pro12}, there is a permutation matrix $P$ such that $P\tilde{B}\tilde{A}P^t=W_2$ or $P\tilde{B}\tilde{A}P^t=J_2$. Let $B=P\tilde{B}$ and $A=\tilde{A}P^t$. Then $AB=\tilde{A}P^tP\tilde{B}=\tilde{A}\tilde{B}=M$.
Thus $A$ and $B$ satisfy condition $\mbox{(i)}$.

Proof of the conditions $\mbox{(ii)}$ and $\mbox{(iii)}$ are similar to the proof of Theorem \ref{theo12}.  $\Box$

By a similar argument, we can reinterpret conditions $\mbox{(i)}$, $\mbox{(ii)}$ and $\mbox{(iii)}$ of Theorem \ref{theo15} to show that if $M$ satisfies $k(M)=2$, then $M$ is one of the $21$ basic types of matrices which we will show in the following.

\begin{theorem} Suppose $M$ is an $n \times n$ Boolean matrix with $b(M)=b=2$. Let $M=AB$ be a Boolean rank factorization. Then $M$ is primitive with $k(M)=2$ if and only if there is a permutation matrix $P$ such that $PMP^t$ has one of the forms in Table $2$ if $BA=W_2$ or $PMP^t$ has one of the forms in Table $3$ if $BA=J_2$. %if $3\leq b\le n-1$ or Table $2$ if $b=2$.
\end{theorem}

  In Table $2$ and Table $3$ %and Table $2$
the rows and columns of each matrix are partitioned conformally, so that each diagonal block is square. % and the top left hand submatrix common to each has $b$ blocks in its partitioning.

$$ Table\  2\ \  (b = 2) $$
\begin{tabular}{lllll} \hline \\
 \ \ \ \ \ \ \ \ \hskip 3cm &  $\left[ \begin{array}{cc|c}
              0 & J & 0     \\
              J & 0 & J \\ \hline
              J & 0 & J
            \end{array} \right],$ &
$\left[ \begin{array}{cc|c}
              0 & J & 0     \\
              J & 0 & J \\ \hline
              J & J & J
            \end{array} \right],$  &
$\left[ \begin{array}{cc|cc}
              0 & J & 0  & 0  \\
              J & 0 & J  & J  \\ \hline
              J & 0 & J  & J \\
              J & J & J  & J
            \end{array} \right].$  & \hskip 3cm \ \  \\ \\
\hline
\end{tabular}

$$ Table\  3\ \  (b = 2) $$
\begin{tabular}{llll} \hline
 $\left[ \begin{array}{cccc}
              J & J & 0 & 0 \\
              0 & 0 & J & J \\
              J & J & 0 & 0 \\
              0 & 0 & J & J
            \end{array} \right],$ &
$\left[ \begin{array}{cccc|c}
              J & J & 0 & 0 & J \\
              0 & 0 & J & J & 0 \\
              J & J & 0 & 0 & J \\
              0 & 0 & J & J & 0 \\ \hline
              J & J & J & J & J
            \end{array} \right],$  &
$\left[ \begin{array}{cccc|c}
              J & J & 0 & 0 & 0 \\
              0 & 0 & J & J & J \\
              J & J & 0 & 0 & 0 \\
              0 & 0 & J & J & J \\ \hline
              J & J & J & J & J
            \end{array} \right],$  &
$\left[ \begin{array}{cccc|cc}
              J & J & 0 & 0 & J & 0 \\
              0 & 0 & J & J & 0 & J \\
              J & J & 0 & 0 & J & 0 \\
              0 & 0 & J & J & 0 & J \\ \hline
              J & J & J & J & J & J \\
              J & J & J & J & J & J
            \end{array} \right],$
\end{tabular}

\begin{tabular}{llll}
 $\left[ \begin{array}{ccc}
              J & J & 0 \\
              0 & 0 & J \\
              J & J & J
            \end{array} \right],$ &
$\left[ \begin{array}{ccc|c}
              J & J & 0 & J \\
              0 & 0 & J & 0 \\
              J & J & J & J \\ \hline
              J & J & J & J
            \end{array} \right],$  &
$\left[ \begin{array}{ccc|c}
              J & J & 0 & 0 \\
              0 & 0 & J & J \\
              J & J & J & J \\ \hline
              J & J & 0 & 0
            \end{array} \right],$  &
$\left[ \begin{array}{ccc|c}
              J & J & 0 & 0 \\
              0 & 0 & J & J \\
              J & J & J & J \\ \hline
              0 & 0 & J & J
            \end{array} \right],$
\end{tabular}

\begin{tabular}{ll}
 $\left[ \begin{array}{ccc|cc}
              J & J & 0 & J & 0 \\
              0 & 0 & J & 0 & J \\
              J & J & J & J & J \\ \hline
              J & J & J & J & J \\
              J & J & 0 & J & 0
             \end{array} \right],$  &
$\left[ \begin{array}{ccc|cc}
              J & J & 0 & J & 0 \\
              0 & 0 & J & 0 & J \\
              J & J & J & J & J \\ \hline
              J & J & J & J & J \\
              0 & 0 & J & 0 & J
             \end{array} \right],$
\end{tabular}

\begin{tabular}{llll}
 $\left[ \begin{array}{ccc}
              J & J & J \\
              J & 0 & 0 \\
              0 & J & J
            \end{array} \right],$ &
$\left[ \begin{array}{ccc|c}
              J & J & J & J \\
              J & 0 & 0 & J \\
              0 & J & J & 0 \\ \hline
              J & 0 & 0 & J
            \end{array} \right],$  &
$\left[ \begin{array}{ccc|c}
              J & J & J & J \\
              J & 0 & 0 & J \\
              0 & J & J & 0 \\ \hline
              0 & J & J & 0
            \end{array} \right],$  &
$\left[ \begin{array}{ccc|c}
              J & J & J & J \\
              J & 0 & 0 & 0 \\
              0 & J & J & J \\ \hline
              J & J & J & J
            \end{array} \right],$
\end{tabular}

\begin{tabular}{llll}
 $\left[ \begin{array}{ccc|cc}
              J & J & J & J & J \\
              J & 0 & 0 & J & 0 \\
              0 & J & J & 0 & J \\ \hline
              J & 0 & 0 & J & 0 \\
              J & J & J & J & J
             \end{array} \right],$  &
$\left[ \begin{array}{ccc|cc}
              J & J & J & J & J \\
              J & 0 & 0 & J & 0 \\
              0 & J & J & 0 & J \\ \hline
              0 & J & J & 0 & J \\
              J & J & J & J & J
             \end{array} \right],$ &
$\left[ \begin{array}{cccc}
              J & J & J & J \\
              J & J & J & J \\
              J & 0 & J & 0 \\
              0 & J & 0 & J
            \end{array} \right],$  &
$\left[ \begin{array}{cccc}
              J & J & J & J \\
              J & J & J & J \\
              0 & J & 0 & J \\
              J & 0 & J & 0
            \end{array} \right].$ \\ \\ \hline
\end{tabular}

\end{document}